\documentclass[12pt]{article} % For LaTeX2e
\usepackage[preprint]{tmlr}

\usepackage{hyperref}
\usepackage{url}

\usepackage{tikz}
\usetikzlibrary{babel}
\usetikzlibrary{patterns}
\usetikzlibrary{decorations,calligraphy}
\usepackage{booktabs}
\usepackage{multirow}

\usepackage{microtype}
\usepackage{graphicx}

\usepackage{hyperref}
\usepackage{subcaption} %  for subfigures environments 

\usepackage{amssymb}
\usepackage{mathtools}
\usepackage{bbm}

\usepackage{algorithm}
\usepackage[noend]{algorithmic}

\usepackage{amsmath,amsfonts}

\title{Decision-Focused Forecasting: \\ A Differentiable \\ Multistage Optimisation Architecture}

\author{\name Egon Peršak \email E.Persak@sms.ed.ac.uk \\
      \addr School of Mathematics\\
      University of Edinburgh
      \aND
      \name Miguel F. Anjos \\
      \addr School of Mathematics\\
      University of Edinburgh
      }

\def\month{MM}  % Insert correct month for camera-ready version
\def\year{YYYY} % Insert correct year for camera-ready version
\def\openreview{\url{https://openreview.net/forum?id=XXXX}} % Insert correct link to OpenReview for camera-ready version

\begin{document}

\maketitle

\begin{abstract}
Most decision-focused learning work has focused on single stage problems whereas many real-world decision problems are more appropriately modelled using multistage optimisation. In multistage problems contextual information is revealed over time, decisions have to be taken sequentially, and decisions now have an intertemporal effect on future decisions. Decision-focused forecasting is a recurrent differentiable optimisation architecture that expresses a fully differentiable multistage optimisation approach. This architecture enables us to account for the intertemporal decision effects of forecasts. We show what gradient adjustments are made to account for the state-path caused by forecasting. We apply the model to multistage problems in energy storage arbitrage and portfolio optimisation and report that our model outperforms existing approaches. 
\end{abstract}

\section{Introduction}
Many real world decision problems are recurrent and decisions taken now have a bearing on future decisions. Since information about such problems is revealed over time it is preferable to model them as having multiple stages so each set of decisions considers the most recently available contextual information. The broad utility of being able to systematically solve such problems is obvious. It is utilised as a core decision making technology in applications as varied as energy systems planning, finance, radiotherapy, and supply chain management. The mathematical optimisation of multistage problems is computationally intense and often relies on specialised decomposition algorithms due to the size of real-world problems and scenario trees needed to capture them. This is further complicated by the fact that problem parameters are often not known before decisions have to be made and that prediction in the time series domain is famously difficult. 

% Alignment problem
The multistage setting inherits the alignment problem of contextual optimisation: the performance on the prediction task is not necessarily aligned with performance on the downstream decision task. Decisions are often highly sensitive to errors in problem parametrisation. This alignment problem is hard to account for a priori to the decision task as we do not have access to the conditional distribution of the unknown parameters. Decision-focused learning emerged to solve the alignment problem problem by training predictive pipelines with respect to downstream decision performance and in doing so implicitly learning the contextual sensitivity structure of the problem. In the case of multistage optimisation the prediction problem is harder and decision taken now propagate to future decision-making exacerbating the alignment problem.

We propose a model that expresses the deterministic policy multistage optimisation approach as a fully differentiable computation to tackle the alignment problem in multistage decision problems. The deterministic policy approach is a simple approach in which only one prediction/forecast is made at each stage. The relative ease of computation in this approach enables the use of differentiable optimisation and consequently an extension of the decision-focused paradigm to multistage problems.

Our main contributions are as follows:
\begin{itemize}
    \item We propose the first decision-focused model for continual multistage optimisation problems. Our model expresses a fully differentiable deterministic policy approach to multistage optimisation.
    \item We analyse the gradients computed by our model to illustrate why this method improves on conventional DFL in multistage settings. In contrast to conventional DFL our approach adjusts for the state-path the prediction model causes enabling intertemporal alignment.
    \item We present two experiments in which we implement our model and show it successfully compares with DFL and the more conventional prediction-focused approach.
\end{itemize}

\section{Background}
\paragraph{Setting}
We focus on structured decision problems under parameter uncertainty which have to be resolved over discrete time. The temporal structure of the process is as follows: at stage $t$ we have an existing state $s_{t-1} \in \mathcal{R}^{|s|}$ and access to context $z_t \in \mathcal{R}^{c_t}$. We need to make decisions $\textbf{x}_t \in \mathcal{R}^d$ before the uncertain parameters  $\theta_t \in \mathcal{R}^p$ for the stage are revealed. This temporal structure is often referred to as here-and-now decisions. We assume the state is certain so not determined by the realisation of $\theta_{t-1}$ and that the uncertainty is only in the objective. The state is a subset of decision variables $s_{t} \subseteq x_{t}$. The context $z_t$ is a catch-all term which includes all past realisations of the unknown parameters and previous auxiliary contexts. The decisions at each stage are constrained to a set $x_t \in \mathcal{X}_t(s_{t-1})$ determined only by the previous state. We restrict ourselves to sets which are expressible with convex optimisation compliant constraints. We assume that this set is non-empty for any value of $s$ also known as relatively complete recourse. The decisions at each stage are scored using a stage-wise objective function determined by the uncertain parameters which we denote as $f(x_t,\theta_t)$. Let the index $t:t+H-1$ denote the matrix of variables corresponding to $H$ stages $t$ to $t+H-1$ for example $x_{t:t+H-1} \in \mathcal{R}^{d \times H}$. Assuming a risk neutral decision maker the decision task for a task horizon of $H$ stages is:
\begin{equation}
    \label{eq: paper multistage}
    \begin{array}{rl}
        \underset{x_{t:t+H-1}}{\min} & \underset{\theta_{t:t+H-1}}{\mathbb{E}} \left[ \sum_{i=0}^{H-1} f(x_{t+i},\theta_{t+i}) \right]  \\
        s.t. & x_{t+i} \in \mathcal{X}_{t+i}(s_{t+i-1}) \quad \forall i \in \{0:H-1\}
    \end{array}
\end{equation}
This problem is difficult to solve even when the exact (conditional) distribution of $\theta_{t:t+h-1}$ is known for all but the simplest distributions. In practice $\theta_{t}$ is a multivariate time series with complex dependency structures making problem (\ref{eq: paper multistage}) intractable. Multistage optimisation is primarily concerned with finding good approximations to this problem. 

\paragraph{Multistage Optimisation}
The multistage optimisation method zoo is vast and there is no dominant method as tasks differ in nature and scope. The expectation is approximated using sampling which can be simple or have more intricate structure such as scenario trees. In multi-scenario approaches a significant feature are non-anticipativity constraints meaning decision variables have to take the same value in each scenario for a given stage. Methods differ by how many scenarios they consider and how they treat non-anticipativity constraints \citep{wang2022impact}. The simplest multistage optimisation approach is having a single scenario at each stage known as the deterministic policy approach. Two-stage approaches relax non-anticipativity constraints for stages after the current stage meaning the decisions taken now reflect flexibility later. Multi-scenario approaches increase the problem size with the number of scenarios in terms on the number of constraints and the number of variables in the case of two-stage approaches. The increased problem size often requires the use of decomposition algorithms. A popular family of approaches constrains possible future decision making to linear decision rules, which results in a significantly easier optimisation problem. 

The distribution of the unknown parameters is typically not known and has to be estimated based on past observations of context and realisations. The forecasting problem with a forecast horizon of $h$ is finding a mapping $m^h(z_t;W) = \hat{\theta}_{t:t+h-1}$ which optimises some loss criterion when compared with true realisations ${\theta}_{t:t+h-1}$. We use $W$ to denote the learnable parameters of the prediction model. This prediction problem is difficult due to the nature of the data generating process (DGP) in time series. The unknown parameters $\theta_t$ are typically highly temporally dependent (either locally or periodically) and the dynamics of the DGP often change over time (non-stationarity). Time series tend to have fewer observations available than other learning tasks. Forecasting is difficult even for single scenario point estimates, forecasting/simulating multiple scenarios is even more challenging.

Given that new contextual information is available at each stage it is sensible to update the forecast and solve the problem approximation at each stage to obtain decisions for the current stage. This also provides a solution to problems which are continual and do not have a finite task horizon. Since solving and infinite horizon problem directly is not computationally tractable we instead tackle it by repeatedly solving the same problem but only for the stages in a set fixed planning horizon. We refer to this clipped version of the problem as the planning problem. We denote the planning horizon as $h$ to distinguish it from the task horizon $h \leq H$. Minor adjustments are often added to the planning problem to reflect the continual nature of the overall task such as forcing the final state to be the same as the initial state or assigning value to the terminal state. Determining the optimal planning horizon in contextual multistage optimisation is a hard problem as it has to balance the myopia of decision planning with the limited predictability of unknown parameters for longer forecasting horizons. 

\paragraph{Decision-Focused Learning}
The key difficulty in contextual optimisation is that the performance on the prediction task does not necessarily align with performance the decision task. This alignment depends on the local sensitivity of the optimal solution to errors in prediction and the unobserved conditional uncertainty of the DGP. Even very low generalisation error in prediction models can result in highly suboptimal decision making. It is hard to account for the alignment factors ex ante analytically due to combinatorially large feasible spaces and lack of information about the DGP. Decision-focused learning (DFL) solves this alignment problem by explicitly training the prediction model for performance on the decision task. DFL models learn the contextual alignment dynamics of the decision problem. We use regret as our measure for decision quality. Let $F$ be the whole objective function and $x^*(\hat{\theta})$ the optimal solution given parameters $\hat{\theta}$, regret is defined as how much objective value is lost as a result of our prediction versus having perfect information:
\begin{equation}
    \mathcal{L}_{regret} =  F(x^*(\hat{\theta}),\theta) - F(x^*(\theta),\theta)
\end{equation}

Mathematically optimising for this loss function is difficult because it requires differentiating across an implicit function $x^*$ which has an $\mathrm{arg}\min$ operator. Differentiable optimisation (DO) has seen a lot of recent attention \citep{amos2017optnet,elmachtoub2022smart, mandi2023decision}. Our examples contend with quadratic programming (QP) problems for which gradients can be obtained analytically by solving a linear system of optimality condition differentials. In our implementation we employ cvxpylayers \citep{cvxpylayers2019}, which extends this method to general conic programs \citep{agrawal2019differentiating} of which QP is a special case. Not all DO methods can be used as a layer as some compute or approximate $\frac{\partial \mathcal{L}}{\partial \hat{\theta}}$ instead of $\frac{\partial x^*}{\partial \hat{\theta}}$ for efficiency or mathematical reasons. We refer to methods which can be used as a layer as implicit layers. \citet{mandi2023decision} present an excellent recent survey on the area. From now on we use the notation $x^*: \mathcal{R}^{p \times h} \mapsto \mathcal{R}^{d \times h}$ for differentiable optimisation layers and $x^*_t \in \mathcal{R}^{d}$ to as the solution for the current state $t$.

In multistage optimisation decisions taken now have an effect on downstream decision making. In our setting this is transmitted by the state variables. The intertemporal effect makes the alignment problem harder as errors in prediction now also have an effect on decisions later. Extending the decision-focused paradigm to multistage problems is therefore highly desirable. This work does not develop new DO methods, instead proposing decision-focused forecasting (DFF) a novel neural architecture composed of many implicit layers which expresses multistage optimisation in a fully differentiable computation. We express the simplest multistage optimisation approach: deterministic policy. This choice is motivated by two factors: the relative difficulty of learning simulators compared to point predictors, and the fact that implicit layers typically require at least one solver call per instance making the larger multi-scenario optimisation problems computationally burdensome. Our architecture is extendable to two-stage multistage approaches with the use of differentiable simulators.

\subsection{Existing Work}
There have been limited attempts to extend DFL to multistage settings. \citet{hu2024two} propose a two-stage approach for solving mixed-integer linear programs with parameter uncertainty in constraints. They focus on problems where recourse is available after the observation of true parameters and model this by adding a second stage penalty function which reflect the cost of the recourse adjustment to feasibility. While originally designed for packing and covering linear programs the framework generalises to the simplest multistage optimisation case of a single scenario with two stages. This approach was later extended in \citet{hu2024multi} to finite multistage settings by training a decision-focused prediction model for each stage. It is not clear how well this non-recursive architecture extends to continual settings.

Deterministic policy approaches are highly reminiscent of model predictive control (MDP) from engineering. Differentiable optimisation has been applied to learn the structure of model predictive control problems \citep{east2020infinite}. The objective of the learning in MDP is a set of parameters directly describing the problem dynamics using imitation learning whereas we aim to learn how to contextually predict problem parameters.

Some multistage stochastic optimisation approaches integrate contextual/covariate/side information. \citet{ban2019dynamic} use a regularised regression model to construct scenario trees informed by contextual information. This approach is not decision-focused but provides great insight in how to construct scenario trees for cold-start forecasting problems. \citet{bertsimas2023dynamic} use non-parametric machine learning approaches to assess similarity between a given context and contexts in the training set and approximates the true problem by weighing past realisations, which are treated as scenarios, based on the similarity score. In addition a norm ball around each scenario is considered to account for aleatoric uncertainty. The theoretical version of this modelling approach has nice asymptotic convergence properties. The problem is approximated using linear decision rules in a two-stage manner allowing for multiple policies. This approach is more decision-focused but completely relies on existing observations being sufficiently similar to future observations. It does not attempt to learn the latent structure of the DGP to generalise to new contexts which is necessary in situations with limited data.

Multistage optimisation problems can be cast as Markov decision problems and therefore solved with (deep) reinforcement learning (RL). Given the structure of the problem is known exactly an RL approach is at a disadvantage since it has to learn the problem dynamics implicitly in value/action/policy functions on top of the alignment dynamics. Generally RL is used in multistage decision problems when the structure of the problem is not known exactly but interactions with the problem environment can be simulated/observed or the problem size is computationally intractable. For a more complete discussion of approaches to sequential decision making we refer the reader to the unified framework proposed by \citet{powell2022unified}.

\section{Decision-Focused Forecasting}
\begin{figure}
    \centering
    \includegraphics[width=0.2\linewidth]{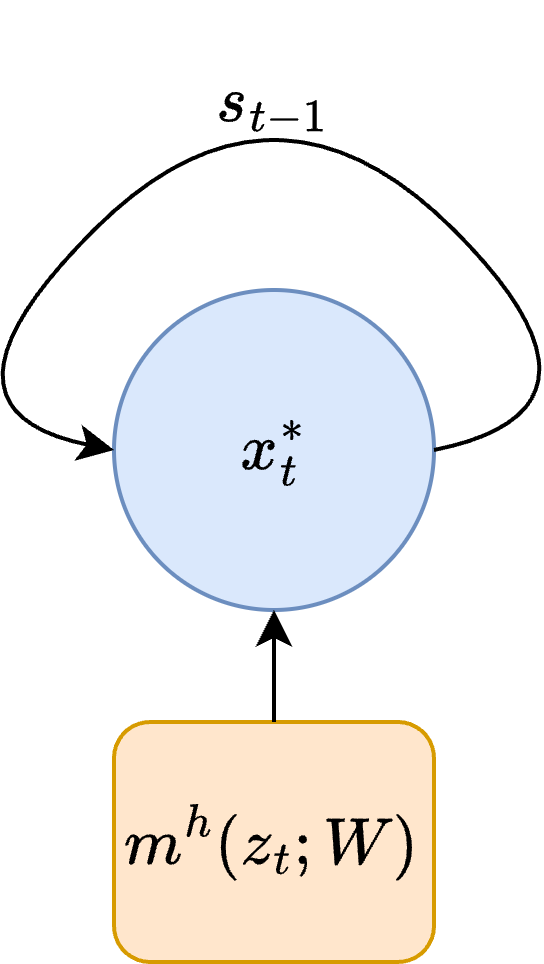}
    \caption{DFF: A deterministic policy approach expressed as a recurrent neural architecture.}
    \label{fig:recurrent DO}
\end{figure}

The iterated optimisation in multistage approaches for continual problems can be expressed as a recurrent neural architecture. In the deterministic policy approach an optimisation problem with a planning horizon of $h$ solved at stage $t$ contains two groups of uncertain parameters: the current state $s_{t-1}$ which is part of the solution to the previous stage and the forecast $\hat{\theta}_{t:t+h-1}$ based on the context $z_t$ at the current stage. We can express this as a neural architecture by using a DO layer as a recurrent layer which takes as inputs a part of the output of the previous layer and the outputs of a prediction model. We name this recurrent architecture decision-focused forecasting (DFF). Note that the DO layer does not contain any learnable parameters and that the learning is restricted to the forecasting model parameterised by $W$. The architecture is illustrated in Figure \ref{fig:recurrent DO}.

We want to derive a policy driven by forecasts which performs well on a continual version of the problem. In most cases we only have one continual observation of the underlying process. The crux of DFF is passing gradient information across stages. The training scheme therefore has to balance having enough different shorter tasks with each task featuring enough stages to express intertemporal effects. We use a sliding window approach to partition a series into $H-h-l+1$ different subtasks. $h$ is the planning horizon and we also use it as our forecasting horizon. The $l$ is the lookback horizon of how many most recent past stages we include as context. For convenience we set the task horizon of each subtask as $h$. Note that the length of each individual subtask horizon can be made longer or shorter than the planning horizon, we did not experiment with what effect that would have on training but speculate that the generalisation of the method could be improved by training on multiple subtask horizons. We fix the initial state as a constant which is the same for all subtasks. This is not wholly representable of the problem and one potential way the training set of subtasks could be bolstered is by generating multiple valid initial states for each subtask. The pseudocode for DFF is presented in algorithm \ref{alg: DFF}. Note the recursive nature of the computation: the decision $x^*_{t+h-1}$ directly depends on the forecast $\hat{\theta}_{t+h-1:t+2h-2}$ and recursively through previous states on all previous forecasts.
% Algorithm
\begin{algorithm}[tb]
   \caption{Epoch of DFF Training}
   \label{alg: DFF}
\begin{algorithmic}
   \STATE {\bfseries Model Components:} forecasting model $m^h(z;W)$, \\ optimisation layer $x^*(s,\hat{\theta})$.
   \STATE {\bfseries Input:} planning horizon $h$, problem horizon $H$, context $z_{t:t+H-1}$, true parameters $\theta_{t:t+H-1}$, batch size $\mathcal{B}$, learning rate $\eta$.
   \REPEAT
   \STATE Draw $\mathcal{B}$ indices from $0$ to $H - 1 - h$ without replacement.
   \STATE Initialise initial states $s_{t+b-1}$ for $b \in \mathcal{B}$.
   \STATE \underline{Forward:} \hfill \COMMENT{Batched for subtasks in $\mathcal{B}$}
   \FOR{$i=0$ {\bfseries to} $h-1$ \hfill}
   % \IF{$x_i > x_{i+1}$}
   \STATE $\psi = t+b+i$.
   \STATE $m^h(z_{\psi}) =\hat{\theta}_{\psi:\psi+h-1}$ \hfill \COMMENT{Forecast}
   \STATE $x^*_{\psi}(s_{\psi-1},\hat{\theta}_{\psi:\psi+h-1}) = x_{\psi}$ \hfill \COMMENT{Optimisation}
   \STATE Save $x_{\psi}$ and $s_{\psi} \subseteq x_{\psi}$.
   % \ENDIF
   \ENDFOR
   \STATE {\bfseries return:} $\textbf{x}_b = x_{t+b:t+b+h-1}$, $\Theta_b = \theta_{t+b:t+b+h-1}$
   \STATE \underline{Backward:}
   \STATE $\mathcal{L} = F(\textbf{x}_b,\Theta_b)$
   \STATE $\omega = t + b + h - 1$
   \STATE $\frac{\partial \mathcal{L}}{\partial W} = 
   \frac{1}{|B|}\sum_{b \in B} \frac{\partial }{\partial W } \left[ 
   \sum_{i=0}^{h-1}
   f(x_{\omega - i},\theta_{w-i})
   \right]
   $
   \STATE where $x_{\omega - i} = x^*_{\omega - i}(s_{\omega - i - 1},m^h(z_{\omega -i},W))$
   \STATE $W \gets W - \eta \frac{\partial \mathcal{L}}{\partial W}$
   \UNTIL{All indices drawn}
\end{algorithmic}
\end{algorithm}

\subsection{Analysis of Gradients}
Differentiating through this sequential, forecast determined optimisation procedure results in much more complex gradients compared to conventional DFL. The prediction model weights are optimised using gradient descent which relies on estimating the gradient $\frac{\partial \mathcal{L}}{\partial W}$. We analyse the computation of this gradient for a subtask starting at stage $t$. Since the objective function is additive we can express it as:
\begin{equation}
\frac{\partial \mathcal{L}}{\partial W} = 
\sum_{i = 0}^{h-1} 
\frac{\partial f(x_{t+i},\theta_{t+i})}{\partial W}
\end{equation}
Denote the forecast for $h$ stages at stage $t+i$ as $\hat{\theta}_{t+i:t+i+h-1} = \hat{\Theta}_{t+i}$. The term for $i=0$ is the most straightforward case as the initial state $s_{t-1}$ is a constant. 
\begin{gather}
        \frac{\partial f(x_{t},\theta_t)}{\partial W} =
\frac{\partial f(x_t^*(m^h(z_t;W),s_{t-1}),\theta_{t})}{\partial W} = \\
= \frac{\partial f(x_{t},\theta_{t})}{\partial x_{t}}
\frac{\partial x_t^*(s_{t-1},\hat{\Theta}_{t})}{\partial W}
\end{gather}
\begin{equation}
\frac{\partial x_t^*}{\partial W} =
\frac{\partial x_t^*(s_{t-1},\hat{\Theta}_{t})}{\partial \hat{\Theta}_{t}}
\frac{\partial m^h(x_t;W)}{\partial W}
\end{equation}
The difficult term to compute is the DO Jacobian $\frac{\partial x^*}{\partial \hat{\Theta}}$: determining how changes in predicted parameters affect the solution to the optimisation problem. Note that only keeping a subset of the decisions from an optimisation layer is equivalent to all other decisions always mapping to $0$ in the loss function and thus not being relevant in the gradient computation. For $i \geq 1$ the input $s_{t+i-1}$ to the optimisation layer recursively depends on previous forecasts and decisions. We state the gradient terms for $i \geq 1$ recursively:
\begin{equation}
\frac{\partial f(x_{t+i},\theta_{t+i})}{\partial W} =
\frac{\partial f(x_{t+i},\theta_{t+i})}{\partial x_{t+i}}
\frac{\partial x^*_{t+i}(s_{t+i-1},\hat{\Theta}_{t+i})}{\partial W} 
\end{equation}
\begin{equation}
\label{eq:recursion}
\frac{\partial x^*_{t+i}}{\partial W} = 
\frac{\partial x^*_{t+i}}{\partial \hat{\Theta}_{t+i}}
\frac{\partial m^h(z_{t+i};W)}{\partial W} + 
\frac{\partial x^*_{t+i}}{\partial s_{t+i-1}}
\frac{\partial s^*_{t+i-1}}{\partial W}
\end{equation}
\paragraph{State-Path Gradient Adjustment} Everything up to equation (\ref{eq:recursion}) is the same as for conventional DFL. Our model differs by adding the term $\frac{\partial x^*_{t+i}}{\partial s_{t+i-1}}
\frac{\partial s^*_{t+i-1}}{\partial W}$. We name this gradient adjustment the state-path adjustment. It accounts for the downstream decision effects of previous forecasting transmitted through the state variables. The state-path adjustment reflects that the forecasts have intertemporal effects. The adjustment term is recursive so the state-path adjustment will consist of an increasing number of terms for later stage decisions in a given subtask. This additional computation is needed to account for the propagation of intertemporal alignment errors in this recursive model. The state-path links the stages and enables the network to learn to induce policies which are future-decision-focused. The adjustment requires additional computation compared to conventional DFL with uncertainty in the objective as the intertemporal gradient $\frac{\partial x^*_{t+i}}{\partial s_{t+i-1}}$ would not be computed.

\paragraph{Intertemporal Gradient} 
For this analysis we restrict ourselves to settings which permit analytic differentiation. Analytic DO gradients are obtained by solving a linear system of the differentials of the KKT optimality conditions. The expensive aspect of this computation is that one solver call is required to obtain the linear system. Solving the linear system to obtain $\frac{\partial x^*_{t+i}}{\partial s_{t+i-1}}$ in addition to the usual DFL gradient $\frac{\partial x^*_{t+i}}{\partial \hat{\Theta}_{t+i}}$ is relatively inexpensive. 

The constraint $x_{t+i} \in \mathcal{X}_{t+i}(s_{t+i-1})$ for the problem at stage $t+i$ is parameterised by the solution to the previous stage. Analytic differentiation of convex optimisation problems with respect to parameters in constraints is well known \citep{amos2017optnet, agrawal2019differentiating}. The constraint set is composed of equalities and inequalities. Inequality constraints can lead to instability in DFF training as they are not necessarily active. In the optimal solution the intertemporal gradient  $\frac{\partial x^*_{t+i}}{\partial s_{t+i-1}^L}$ for the subset of state variables $s_{t+i-1}^L \subseteq s_{t+i-1}$ which are not part of any active constraints is a zero matrix. This is because small perturbations to those constraints will not change the optimal solution. Conversely, the intertemporal gradient is almost always a non-zero matrix when the constraints are active. If the optimal solution given true parameters has a non-empty $s_{t+i-1}^L$ then there is no intertemporal effect for $s_{t+i-1}^L$ at that stage. However, if the predicted parameters result in a solution in which some of $s_{t+i-1}^L$ are part of active constraints the estimated intertemporal gradient for those variables is non-zero. This mismatch is likely given a random initialisation of the prediction model. The recursive nature of the state-path adjustment means the effect of this mismatch is potentially substantial and can lead to instability in training. We do not implement this but note that at the cost of one extra solver call per stage (substantial) the intertemporal gradient at that stage could be set to $0$ for state variables not exhibiting an intertemporal effect in the optimal solution given true parameters.

\section{Experiments}
% https://github.com/EgoPer/Decision-Focused-Forecasting
% https://anonymous.4open.science/r/Decision-Focused-Forecasting-9BCD
\subsection{Battery Storage}
The battery storage arbitrage task is an evolved version of the third experiment in \citet{donti2017task}. The code implementation is available \href{https://github.com/EgoPer/Decision-Focused-Forecasting}{here}. At each stage a decision has to be made about how much to charge or discharge a battery over the next hour in advance of knowing the market price. The price needs to be predicted based on the available context at each stage. The objective is to maximise trading profit over time with problem specific penalties. The original experiment applied DFL in the setting of a single stage stochastic programming approach with an assumption on the distribution of prices. Our version evolves the experiment to a continual setting by adding a term which approximates terminal value at the end of the planning horizon and evaluates the resulting policy on a rolling basis over the test set.
\paragraph{Task Details}
We set a planning horizon of $h=48$. At each stage $t$ we need to decide how much energy to charge $\textbf{x}_t^{in}$ or discharge $\textbf{x}_t^{out}$. This determines our state $\textbf{s}_t$ which reflects how much energy is stored in the battery at the end of the stage. In addition to maximising profit the objective contains quadratic penalties on the size of $\textbf{x}_t^{in},\textbf{x}_t^{out}$ as larger movements are worse for battery health, a quadratic penalty on how far the state is from being half-full which acts as guidance to promote being in more reactive states, and an adjustment to reflect the continual nature of the problem. The planning problem $x^*$ at stage $t$ given predicted parameters $\hat{\theta}_{t:t+h-1}$ and the initial state $s_{t-1}$ is:

\begin{equation}
\begin{array}{cl}
\displaystyle 
x^*(s_{t-1},\hat{\Theta}_t) = 
\mathrm{arg}\min_{\textbf{s},\textbf{x}^{in},\textbf{x}^{out}}  & \sum_{i=0}^{h-1} ( \hat{\theta}_{t+i} [\textbf{x}_{i}^{in} - \textbf{x}_{i}^{out}] \\[4pt]
& + \lambda[\textbf{s}_i-\frac{B}{2}]^2 + \epsilon [\textbf{x}_{i}^{in}]^2 \\[4pt]
& +\epsilon [\textbf{x}_{i}^{out}]^2
) \\[4pt]
& + \hat{\theta}_{t+h-1} \textbf{s}_{t+h-1}
\\[4pt]
s.t. & \textbf{s}_{-1} = s_{t-1} \\[4pt]
\forall i \in 0,...,h-1 
& \textbf{s}_i = \textbf{s}_{i-1} + \gamma \textbf{x}_{i}^{in} - \textbf{x}_{i}^{out}\\[4pt]
& 0 \leq \textbf{s}_i \leq B \\[4pt] 
& 0 \leq z_{i}^{in} \leq c^{in},\\[4pt] 
& 0 \leq z_{i}^{out} \leq c^{out} \: .
\end{array}
\end{equation}
Where $\lambda = 0.1$ is the penalty parameter for being away from a flexible state, $\epsilon = 0.05$ penalises faster charging and discharging in the name of battery health, and $\gamma = 0.9$ is the charging efficiency parameter. The final term in the objective is the value of the energy currently stored at market prices, an approximation of the value of the terminal state. The equality constraints describe the charging and discharging dynamics from stage to stage. The inequality constraints describe the limits to charging $c^{in} = 0.5$, discharging $c^{out} = 0.2$ and the total battery capacity $B = 1$. We set the initial state value of each planning problem as $\frac{B}{2}$.

We predict the prices $\hat{\Theta}_{t} \in \mathcal{R}^{48}$ based on a contextual input $z_t$ which is composed of observed past and some predicted variables of prices and temperature. Another important hyperparameter is the lookback horizon $l$: the number of most recent stages from which parameter realisations are taken as context. In this task we set $l$ to be the same as the forecast horizon $h$. All inputs are normalised before being passed into the neural network $m^h(z_t; W) = \hat{\Theta}_t$. We use a simple feedforward neural network with two hidden layers. We set the latent dimension to be $2h$, use ReLU activation functions, and dropout of $0.1$ after the hidden layers. This is by no means an endorsement of such networks for forecasting in practice. DFF can be employed with any differentiable forecasting architecture.

\paragraph{Training and Evaluation}
From the 6 years of available hourly data for electricity prices we use the first $0.8$ of the data for training and the remainder for testing leading to approximately $42k$ subtasks in the training set. While the training is performed on subtasks, the evaluation on the test set is for the policy produced by the model across the entire task horizon. The test problem contains just over $10k$ stages or just over 14 months of energy storage arbitrage.

\begin{figure}
    \centering
    \includegraphics[width=1\linewidth]{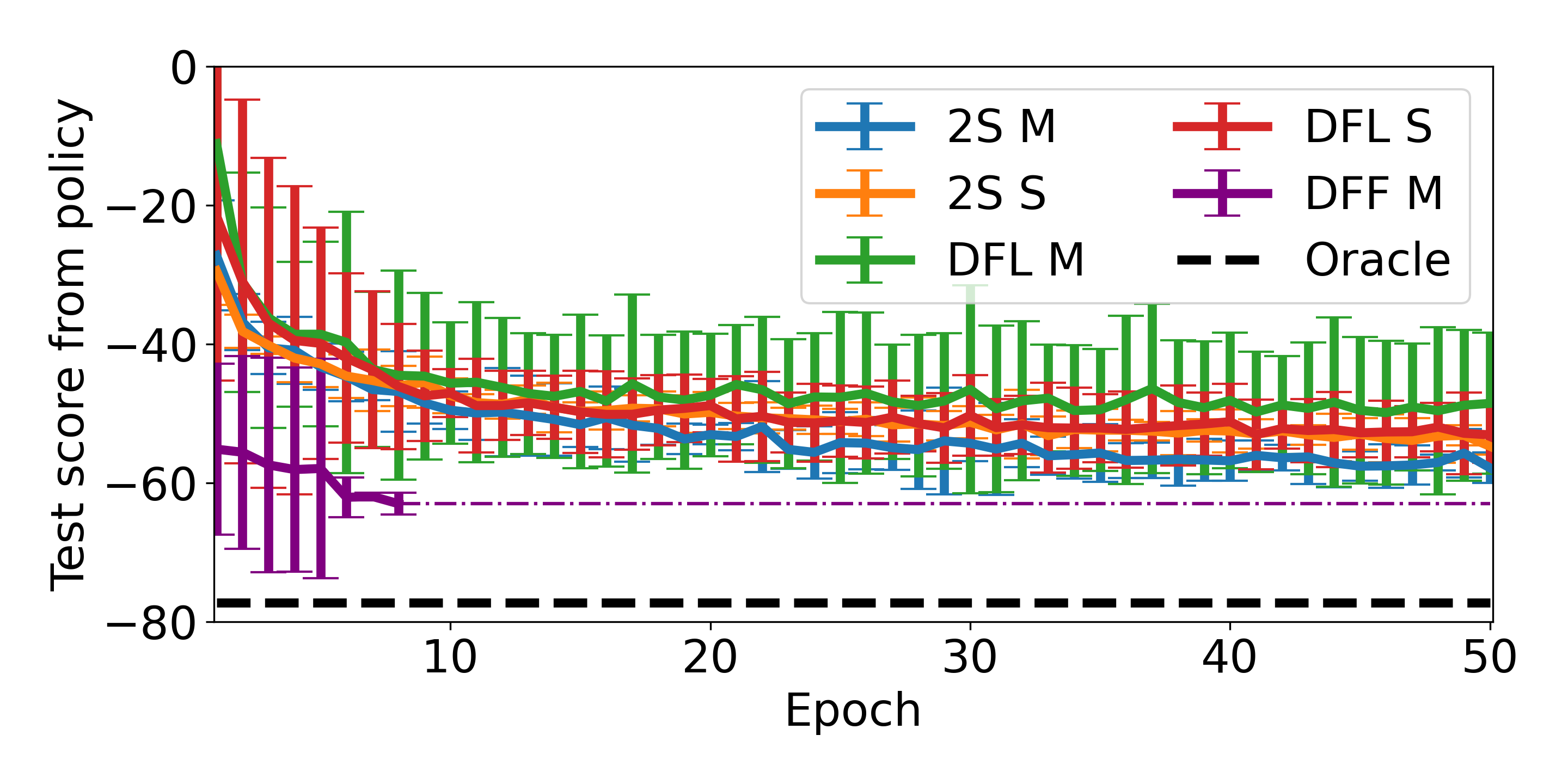}
    \caption{Development of test losses over training on the battery storage arbitrage task. 2-sigma error bars from five training runs. 2S denotes a two stage approach. The M refers to multistage optimisation, the S refers to single-stage optimisation. DFF was trained for 8 epochs, the other two for 50 to balance the relative intensity of computation.}
    \label{fig:results_bs}
\end{figure}

For DFF we train the model for $8$ epochs on an Apple M2 using the Adam optimiser with a batch size of $64$ subtasks using a learning rate of $1\rm{e-}3$. In addition we train two more prediction models with the same neural architecture: a two-stage (2S) model which is optimised for RMSE scored prediction performance, and a DFL model optimised for empirical regret but treating each subtask as a separate single stage optimisation problem. We evaluate the policies from the models for comparison in two ways: 'S' means we only make decisions every $h$ stages effectively partitioning the continual problem into separate single stage problems and 'M' reflecting the same multistage approach employed in DFF where only the decisions for a given stage are made at that stage. This is to check that there is no bias against DFL as it was not trained for the M setting. A training run of $8$ epochs across all three models took about a day on average. An epoch of DFF training took approximately five times as long as DFL training, whereas the training time for 2S is insignificant in comparison. To account for the additional computation performed in DFF due to its recursive gradients we trained 2S and DFL for a total of $50$ epochs.
\paragraph{Results}
We present the results for the task in Figure \ref{fig:results_bs}. We performed the experiment 5 times and display the 2-sigma error bars. DFF trains much faster. This can be explained by the state-path adjustments which can be interpreted as increasing the batch size. Even with the additional $42$ epochs the other two approaches do not even begin to approach the DFF results. The variance of DFF policy results diminishes significantly after about five epochs potentially indicating some form of convergence on this task. 2S approaches outperform DFL approaches which indicates that intertemporal alignment is important for this task. In this experiment mean DFF relative regret is $6.7 \%$ lower than the best performing alternative 2S M. The performance gains from training would likely continue across all models, but with decreasing marginal gains.

Note that the battery storage arbitrage problem is an ideal problem for DFF: it exhibits a high degree of intertemporality, its core dynamics are governed by equality constraints, and the stochastic process of the market clearing price for electricity has a relatively high degree of predictability. It is also made simpler by only having one uncertain parameter in every stage.

\begin{figure}
    \centering
    \includegraphics[width=1\linewidth]{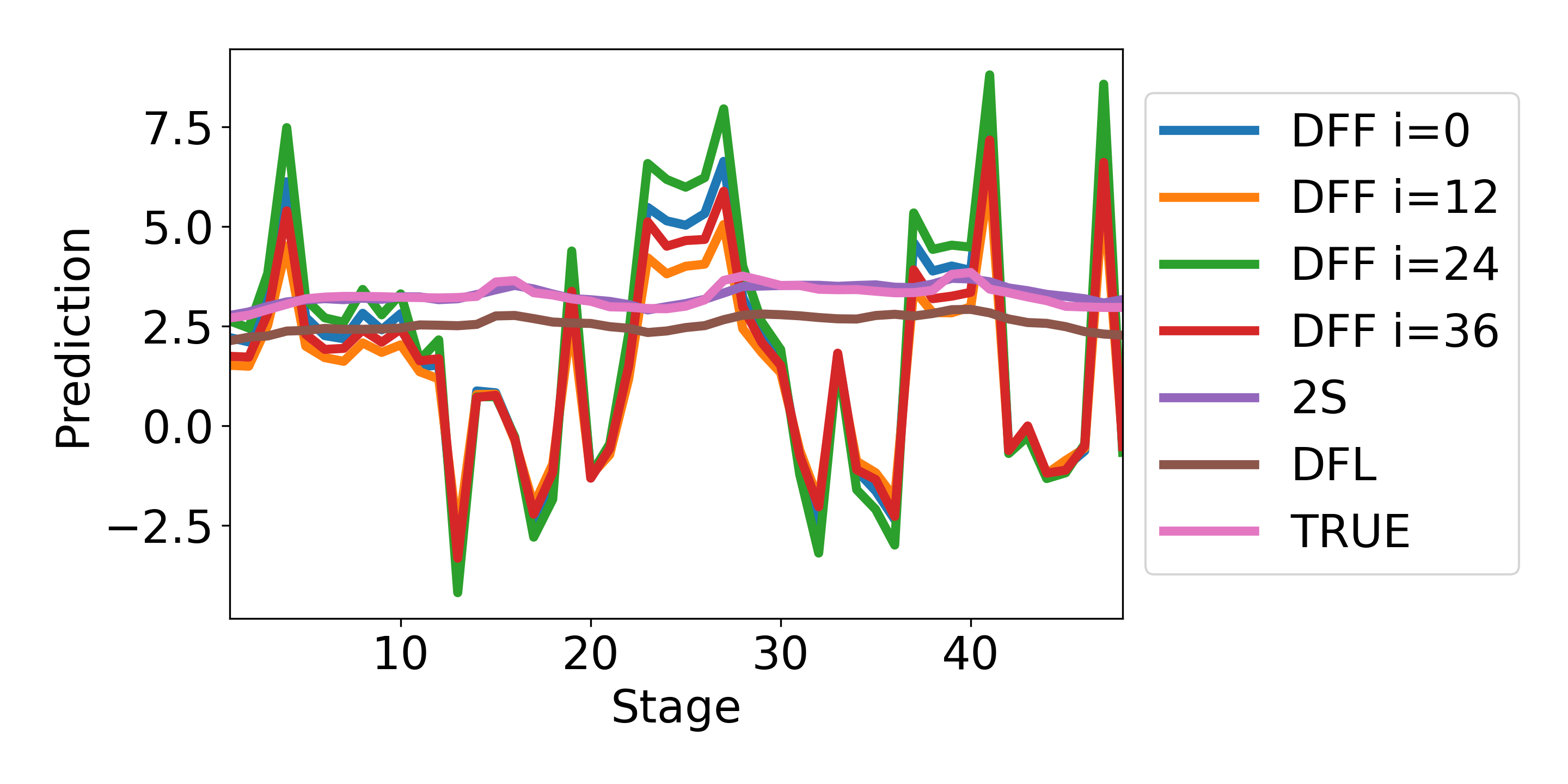}
    \caption{Forecasts produced by the tested model for the battery storage task following one of the experiments. For DFF we also plot the forecasts for contexts $i$ stages after the original one to demonstrate the relative stability of the representation, a phenomenon not shared with other models.}
    \label{fig:forecasts_bs}
\end{figure}

\paragraph{Analysis of Forecasts}
We present the forecasts produced by one of the training runs in Figure \ref{fig:forecasts_bs}. The salient feature is that the forecasts produced by DFF are completely alien in scale and shape. DFL while somewhat off in scale still largely mimics the shape of the 2S and true values. The forecasts produced by DFF appear to mostly preserve their shape over stages. This relative stability phenomenon persists across all five experiments with each exhibiting a unique shape.

Taking the perspective that optimisation is a layer in a broader neural architecture mapping context to decisions the forecasts can be interpreted as latent representations. The representations that decision-focused approaches create to parametrise decision problems reflect something about the structure of the contextual optimisation problem. \citet{mandi2022decision} observe that DFL applied to combinatorial optimisation can be equivalently cast as a learning to rank problem. The underlying structure expressed by decision-focused representations is not easily interpretable, especially for complex settings like multistage optimisation. Our hypothesis is that on this task DFF finds a stable guiding representation which is possible due to the relative predictability of the DGP or the relative simplicity of the optimisation problem. The model effectively finds a high performing local minimum on the planning problem. Exploring what patterns emerge for different multistage problems may produce new understanding about their underlying structure.

\subsection{Portfolio Optimisation}
Our second experiment is a portfolio optimisation problem. In contrast with the battery storage problem, this problem is not ideal for DFF: it has a relatively low degree of predictability, and its intertemporality is transient. We use a slightly modified version of the dataset from \citet{hoseinzade2019cnnpred} which compiles the daily returns from five large stock indices along with a range of contextual data in the form of economic and technical indices. We set a planning horizon of $h = 7$ days and a lookback horizon of $l = 21$ days. The decision task is to allocate portfolio weights $\textbf{s}_t \in \mathcal{R}^5$ across the five indices each day such that the sum of returns over time is maximised less a variance penalty given a set of trading constraints. The planning problem is at stage $t$ is:
\begin{equation}
\label{eq: PO}
\begin{array}{cl}
\displaystyle 
x^*(s_{t-1},\hat{\Theta}_t) = 
\mathrm{arg}\min_{\textbf{s},\Delta}  & \sum_{i=0}^{h-1} -\hat{\theta}_{t+i} \textbf{s}_{i}  + \gamma \textbf{s}_{i}^T \Sigma \textbf{s}_{i} \\
s.t. & \textbf{s}_{-1}= s_{t-1} \\[4pt]
\forall i \in 0,...,h-1 
& \Delta_i = \textbf{s}_{i} - \textbf{s}_{i-1}\\[4pt]
& \mathbbm{1} \textbf{s}_{i} \leq 1 \quad \textbf{s}_{i} \geq 0 \\[4pt]
& \| \Delta_i \|_2 \leq \Lambda \: .
\end{array}
\end{equation}
The covariance matrix $\Sigma$ is fixed and calculated from the training set. We set the variance penalty $\gamma= \frac{1}{2}$. The variable $\Delta_i$ tracks the asset-wise changes from stage to stage. The trading constraints in the penultimate line prohibit shorting or leverage. The final trading constraint is that the change in portfolio from period to period can at most have a euclidean distance of $\Lambda = 0.2$. This is a highly simplified and relatively arbitrary restriction on large movements in the composition in the portfolio. Without this constraint this problem is not intertemporal as nothing else links the stages. We note that it is an inequality constraint which may lead to some instability in the training as outlined in the previous section. Each subtask is initialised with a portfolio of equal weights.

\begin{figure}
    \centering
    \includegraphics[width=1\linewidth]{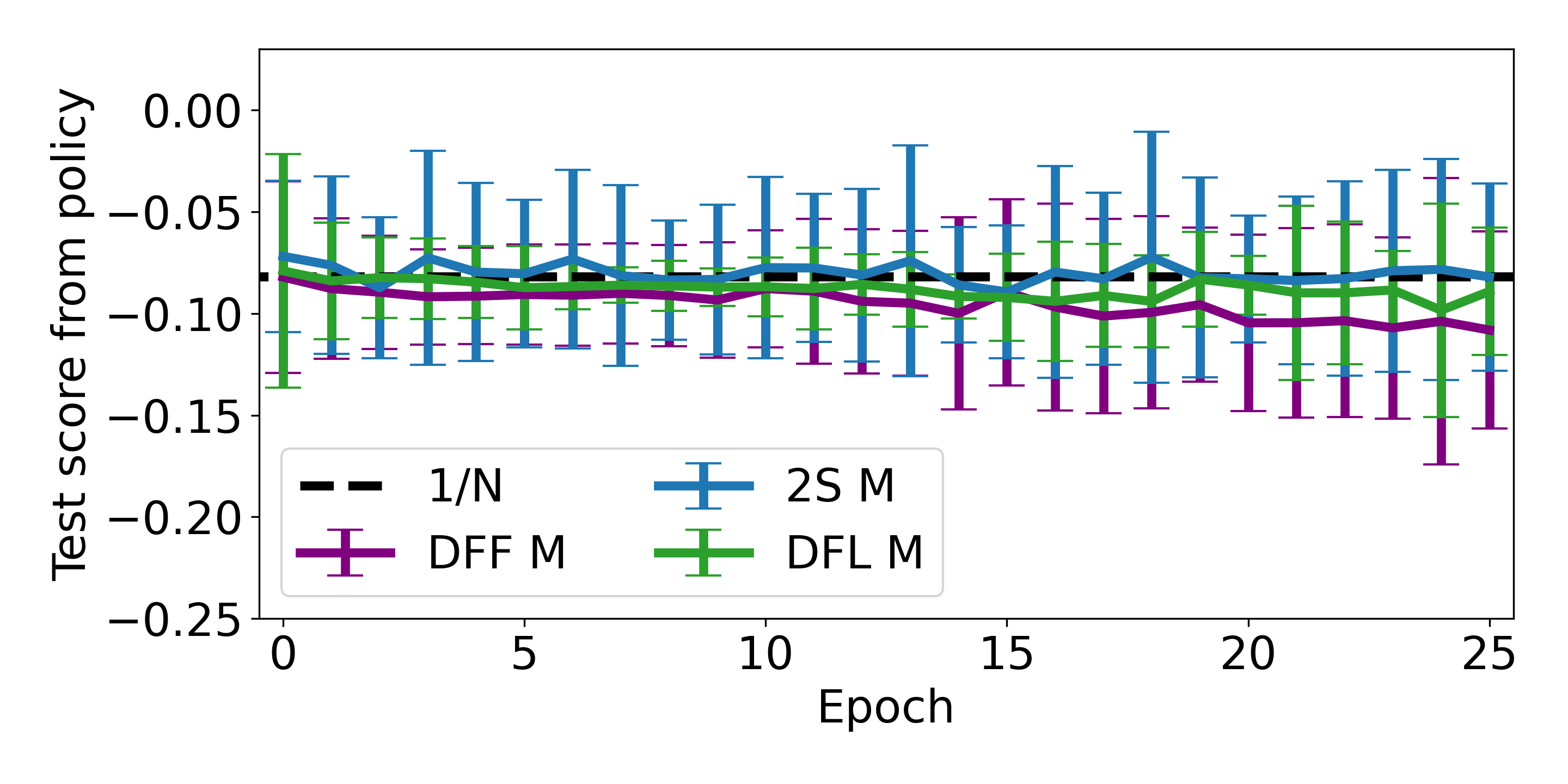}
    \caption{Development of test losses over training on the portfolio optimisation task. 2-sigma error bars from $15$ training runs. $1/N$ denotes the performance of a portfolio with equal weights.}
    \label{fig:results_po}
\end{figure}

In this task we predict the returns of the $5$ individual assets for the next $7$ days so $\hat{\Theta}_t \in \mathcal{R}^{7 \times 5}$. We use a more complex neural network with eight hidden layers which alternate between dropout and layer normalisation. The implicit regularisation is added to stabilise training as the DGP for asset returns exhibits a high degree of noise.
\paragraph{Training and Evaluation}
The training set contains $1.4k$ subtasks which are evaluated on a test task of the final 126 days approximately corresponding to half of the trading days in a year. We use a learning rate of $1e-4$ for the 2S and DFL training, and divide the learning rate by $\sqrt{h}$ for DFF to roughly account for the increase in the size of the gradients as a result of the state-path adjustment and stabilise training. We train all models for $25$ epochs over $15$ training runs. In contrast with the battery storage task the improvements in 2S and DFL performance do not seem to occur at a different pace to DFF which may be explainable by the learning rate adjustment.
\paragraph{Results}
We present the evolution of performance over training in Figure \ref{fig:results_po}. We note the objective value from simply keeping the original equally weighted portfolio (EWP). The mean test task scores from 2S oscillate around the EWP result, not exhibiting a meaningful improvement in performance. DFL shows improvement more often and some individual results are competitive. In contrast, mean DFF performance steadily improves after the tenth epoch. The performance of DFF on this task is evidence that the approach has merit on tasks with less friendly DGPs. However, the variance of DFF results is high and more work is needed to stabilise training.

The relative stability phenomenon between DFF forecasts as observed in the battery storage task appears to a lesser extent in this task: rough high level patterns can be deduced. This can be interpreted to reflect the relative difficulty of this task and the relative importance of contextual information.

\section{Discussion}

% training stability
\paragraph{Training Stability} The training stability in the portfolio optimisation task is poor. This could be due to the inequality constraint defining the intertemporality which can be tackled by looking at whether the constraint is binding in the optimal solution given true parameters and updating the state-path gradient adjustment accordingly. 

The choice of setting the subtask horizon to equal the planning horizon is likely suboptimal. Training could be enhanced by including subtasks of varying lengths in the training set. Furthermore, initialising subtasks with different feasible initial states should improve training by increasing the diversity of subtasks in the training set. These two directions alone can increase the size of a training set by several orders of magnitude.

A natural algorithmic development would be to extend the decision-focused paradigm to multi-scenario multistage optimisation. This would enable greater representational power on the predictive side, enable more hedging in the planning problem, but make each planning problem substantially more expensive computationally. One potential implementation is to generate forecasts/scenarios with different planning horizons and solve them in the same planning problem enabling predictive models to specialise for different degrees of foresight.

\paragraph{Computational Cost} The model we propose is significantly slower in training and has unfavourable scaling properties in terms of the size of the planning problem. Note that this is not a problem at test time. Compared to DFL (which is already slow compared to 2S), each subtask requires sequentially solving $h$ planning problems. For the battery storage problem, we solve approximately $2$ million planning problems each epoch. Implicit layers cannot be parallelised in the same way as ordinary layers due to the sequential nature of mathematical optimisation. On the other hand, the planning problems solved in DFF have a high degree of similarity. Amortised optimisation methods \citep{amos2023tutorial} aim to infer the optimal solution of a problem from parameter inputs based on a set of problem solution observations presenting a potential way of reducing computational expense in DFF training.

% \paragraph{Intertemporality} The relative performance of DFL and DFF can serve as an empirical test for intertemporality in the optimal solution.

\section{Conclusion}
We propose and positively evaluate the first decision-focused approach for continual multistage optimisation. The recurrent implicit layer architecture enables a prediction model to account for the intertemporal effect of forecasts made now on decisions made later. While DFF training is computationally expensive early results suggest a significant improvement in performance on certain continual decision problems. The key areas for improvement are algorithmic innovations for training stability, extending the decision-focused paradigm to more multistage optimisation approaches, and less computationally intense DO techniques. Decision-focused forecasting produces representations which reveal a fascinating underlying structure of contextual multistage optimisation, the precise nature of which is yet to be determined.

\bibliography{bibliography}

\begin{thebibliography}{16}
\providecommand{\natexlab}[1]{#1}
\providecommand{\url}[1]{\texttt{#1}}
\expandafter\ifx\csname urlstyle\endcsname\relax
  \providecommand{\doi}[1]{doi: #1}\else
  \providecommand{\doi}{doi: \begingroup \urlstyle{rm}\Url}\fi

\bibitem[Agrawal et~al.(2019{\natexlab{a}})Agrawal, Amos, Barratt, Boyd, Diamond, and Kolter]{cvxpylayers2019}
A.~Agrawal, B.~Amos, S.~Barratt, S.~Boyd, S.~Diamond, and Z.~Kolter.
\newblock Differentiable convex optimization layers.
\newblock In \emph{Advances in Neural Information Processing Systems}, 2019{\natexlab{a}}.

\bibitem[Agrawal et~al.(2019{\natexlab{b}})Agrawal, Barratt, Boyd, Busseti, and Moursi]{agrawal2019differentiating}
Akshay Agrawal, Shane Barratt, Stephen Boyd, Enzo Busseti, and Walaa~M Moursi.
\newblock Differentiating through a cone program.
\newblock \emph{arXiv preprint arXiv:1904.09043}, 2019{\natexlab{b}}.

\bibitem[Amos \& Kolter(2017)Amos and Kolter]{amos2017optnet}
Brandon Amos and J~Zico Kolter.
\newblock Optnet: Differentiable optimization as a layer in neural networks.
\newblock In \emph{International Conference on Machine Learning}, pp.\  136--145. PMLR, 2017.

\bibitem[Amos et~al.(2023)]{amos2023tutorial}
Brandon Amos et~al.
\newblock Tutorial on amortized optimization.
\newblock \emph{Foundations and Trends{\textregistered} in Machine Learning}, 16\penalty0 (5):\penalty0 592--732, 2023.

\bibitem[Ban et~al.(2019)Ban, Gallien, and Mersereau]{ban2019dynamic}
Gah-Yi Ban, J{\'e}r{\'e}mie Gallien, and Adam~J Mersereau.
\newblock Dynamic procurement of new products with covariate information: The residual tree method.
\newblock \emph{Manufacturing \& Service Operations Management}, 21\penalty0 (4):\penalty0 798--815, 2019.

\bibitem[Bertsimas et~al.(2023)Bertsimas, McCord, and Sturt]{bertsimas2023dynamic}
Dimitris Bertsimas, Christopher McCord, and Bradley Sturt.
\newblock Dynamic optimization with side information.
\newblock \emph{European Journal of Operational Research}, 304\penalty0 (2):\penalty0 634--651, 2023.

\bibitem[Donti et~al.(2017)Donti, Amos, and Kolter]{donti2017task}
Priya Donti, Brandon Amos, and J~Zico Kolter.
\newblock Task-based end-to-end model learning in stochastic optimization.
\newblock \emph{Advances in neural information processing systems}, 30, 2017.

\bibitem[East et~al.(2020)East, Gallieri, Masci, Koutnik, and Cannon]{east2020infinite}
Sebastian East, Marco Gallieri, Jonathan Masci, Jan Koutnik, and Mark Cannon.
\newblock Infinite-horizon differentiable model predictive control.
\newblock \emph{arXiv preprint arXiv:2001.02244}, 2020.

\bibitem[Elmachtoub \& Grigas(2022)Elmachtoub and Grigas]{elmachtoub2022smart}
Adam~N Elmachtoub and Paul Grigas.
\newblock Smart “predict, then optimize”.
\newblock \emph{Management Science}, 68\penalty0 (1):\penalty0 9--26, 2022.

\bibitem[Hoseinzade \& Haratizadeh(2019)Hoseinzade and Haratizadeh]{hoseinzade2019cnnpred}
Ehsan Hoseinzade and Saman Haratizadeh.
\newblock Cnnpred: Cnn-based stock market prediction using a diverse set of variables.
\newblock \emph{Expert Systems with Applications}, 129:\penalty0 273--285, 2019.

\bibitem[Hu et~al.(2024{\natexlab{a}})Hu, Lee, and Lee]{hu2024two}
Xinyi Hu, Jasper Lee, and Jimmy Lee.
\newblock Two-stage predict+ optimize for milps with unknown parameters in constraints.
\newblock \emph{Advances in Neural Information Processing Systems}, 36, 2024{\natexlab{a}}.

\bibitem[Hu et~al.(2024{\natexlab{b}})Hu, Lee, Lee, and Stuckey]{hu2024multi}
Xinyi Hu, Jasper Lee, Jimmy Lee, and Peter Stuckey.
\newblock Multi-stage predict+ optimize for (mixed integer) linear programs.
\newblock \emph{Advances in Neural Information Processing Systems}, 37:\penalty0 64794--64827, 2024{\natexlab{b}}.

\bibitem[Mandi et~al.(2022)Mandi, Bucarey, Tchomba, and Guns]{mandi2022decision}
Jayanta Mandi, V{\i}ctor Bucarey, Maxime Mulamba~Ke Tchomba, and Tias Guns.
\newblock Decision-focused learning: Through the lens of learning to rank.
\newblock In \emph{International conference on machine learning}, pp.\  14935--14947. PMLR, 2022.

\bibitem[Mandi et~al.(2023)Mandi, Kotary, Berden, Mulamba, Bucarey, Guns, and Fioretto]{mandi2023decision}
Jayanta Mandi, James Kotary, Senne Berden, Maxime Mulamba, Victor Bucarey, Tias Guns, and Ferdinando Fioretto.
\newblock Decision-focused learning: Foundations, state of the art, benchmark and future opportunities.
\newblock \emph{arXiv preprint arXiv:2307.13565}, 2023.

\bibitem[Powell(2022)]{powell2022unified}
Warren~B Powell.
\newblock \emph{Reinforcement Learning and Stochastic Optimization: A Unified Framework for Sequential Decisions}.
\newblock Wiley, 2022.
\newblock ISBN 9781119815051.

\bibitem[Wang et~al.(2022)Wang, Cevik, and Bodur]{wang2022impact}
Juyoung Wang, Mucahit Cevik, and Merve Bodur.
\newblock On the impact of deep learning-based time-series forecasts on multistage stochastic programming policies.
\newblock \emph{INFOR: Information Systems and Operational Research}, 60\penalty0 (2):\penalty0 133--164, 2022.

\end{thebibliography}
\bibliographystyle{tmlr}

%%%%%%%%%%%%%%%%%%%%%%%%%%%%%%%%%%%%%%%%%%%%%%%%%%%%%%%%%%%%%%%%%%%%%%%%%%%%%%%
%%%%%%%%%%%%%%%%%%%%%%%%%%%%%%%%%%%%%%%%%%%%%%%%%%%%%%%%%%%%%%%%%%%%%%%%%%%%%%%
% APPENDIX
%%%%%%%%%%%%%%%%%%%%%%%%%%%%%%%%%%%%%%%%%%%%%%%%%%%%%%%%%%%%%%%%%%%%%%%%%%%%%%%
%%%%%%%%%%%%%%%%%%%%%%%%%%%%%%%%%%%%%%%%%%%%%%%%%%%%%%%%%%%%%%%%%%%%%%%%%%%%%%%
% \newpage
% \appendix
% \onecolumn

% You can have as much text here as you want. The main body must be at most $8$ pages long.
% For the final version, one more page can be added.
% If you want, you can use an appendix like this one.  

% The $\mathtt{\backslash onecolumn}$ command above can be kept in place if you prefer a one-column appendix, or can be removed if you prefer a two-column appendix.  Apart from this possible change, the style (font size, spacing, margins, page numbering, etc.) should be kept the same as the main body.
%%%%%%%%%%%%%%%%%%%%%%%%%%%%%%%%%%%%%%%%%%%%%%%%%%%%%%%%%%%%%%%%%%%%%%%%%%%%%%%
%%%%%%%%%%%%%%%%%%%%%%%%%%%%%%%%%%%%%%%%%%%%%%%%%%%%%%%%%%%%%%%%%%%%%%%%%%%%%%%

\end{document}